\documentclass{article}
\usepackage{latexsym}
\usepackage{amsmath}
\usepackage{verbatim}
\usepackage{graphics}
\usepackage{amsfonts}
\usepackage{graphicx}
\usepackage{epstopdf}
\newtheorem{thm}{Theorem}

\newtheorem{prop}{Proposition}
\newtheorem{lem}{Lemma}
\newtheorem{defn}{Definition}
\newtheorem{cor}{Corollary}
\newtheorem{quest}{Question}

\begin{document}

\title{Chord Diagrams and Gauss Codes for Graphs}
\author{Thomas Fleming\\
	Department of Mathematics\\
   	University of California, San Diego\\
   	La Jolla, CA 92093-0112\\
   	{\it  tfleming@math.ucsd.edu}
   \and
   	Blake Mellor\\
	Mathematics Department\\
	Loyola Marymount University\\
	Los Angeles, CA  90045-2659\\
   	{\it  bmellor@lmu.edu}}
\date{}
\maketitle

\begin{abstract}
Chord diagrams on circles and their intersection graphs (also known as {\it circle graphs}) have been intensively studied, and have many applications to the study of knots and knot invariants, among others.  However, chord diagrams on more general graphs have not been studied, and are potentially equally valuable in the study of spatial graphs.  We will define chord diagrams for planar embeddings of planar graphs and their intersection graphs, and prove some basic results.  Then, as an application, we will introduce Gauss codes for immersions of graphs in the plane and give algorithms to determine whether a particular crossing sequence is realizable as the Gauss code of an immersed graph.

\end{abstract}
\tableofcontents

\section{Introduction} \label{S:intro}

Classically, a {\it chord diagram} is a collection of chords of a circle; the intersection graph for these chords is called a {\it circle graph}.  Circle graphs have been intensively studied in graph theory, with classifications given by Even and Itai \cite{ei} and Bouchet \cite{bo}, among others.  In recent years, interest in chord diagrams has spread to topologists as part of the theory of finite type knot and link invariants \cite{bn, bl, cdl}, as well as the theory of virtual knots \cite{ka}.  A natural extension of knot theory is to look at spatial embeddings of more complex graphs, and it is natural to ask whether some analogue of chord diagrams could be equally useful in this context.  The goal of this paper is to define a reasonable notion of chord diagram for general graphs, and to construct some tools to study them.  In particular, we will look at a variety of intersection graphs for these chord diagrams.  We will use these intersection graphs to determine when an embedding of a planar graph can be extended to an embedding of a chord diagram on that graph.  As an application, we will define {\it Gauss codes} for immersions of arbitrary graphs and give algorithms for determining whether a crossing sequence is realizable as a Gauss code.

{\sc Acknowledgement:}  The authors would like to acknowledge the hospitality of Waseda University, Tokyo,
and Professor Kouki Taniyama during the International Workshop on Knots and Links in a Spatial Graph in
July, 2004, where the idea for this project was conceived.  The second author was supported by an LMU Faculty Research grant.

\section{Chord Diagrams for Graphs} \label{S:chord}

\subsection{Definition of chord diagrams} \label{SS:defchord}

A chord diagram for a circle may be characterized as simply a circle with a set of labeled points indicating the endpoints of the chords.  This is the idea we use to define a chord diagram for a general graph.

\begin{defn} \label{D:chord}
Let G be a graph with edge set E and vertex set V.  A {\bf chord diagram on G of degree n} is a collection of 2n points in $G - V$ (so each point is in the interior of an edge), each labeled from a set $\{c_1,...,c_n\}$ such that each label is used exactly twice.  The two points labeled $c_i$ are called the {\bf endpoints} of chord $c_i$.
\end{defn}

Often, we will want to look at {\it oriented} chord diagrams of directed graphs.  This allows us to record whether the chord is on the same or opposite sides of the edges at each endpoint.

\begin{defn} \label{D:orientchord}
Let G be a {\bf directed} graph with edge set E and vertex set V.  An {\bf oriented chord diagram on G of degree n} is a collection of 2n points in $G - V$ , each labeled from a set $\{c_1, \overline{c_1}, c_2, \overline{c_2},..., c_n, \overline{c_n}\}$ such that exactly two points are labeled from each set $\{c_i, \overline{c_i}\}$.  These two points are called the {\bf endpoints} of chord $c_i$; if they have the same label, the endpoints are said to have the same orientation.  If one is labeled $c_i$ and the other $\overline{c_i}$, they are said to have opposite orientations.
\end{defn}

\begin{figure} [h]
    $$\includegraphics{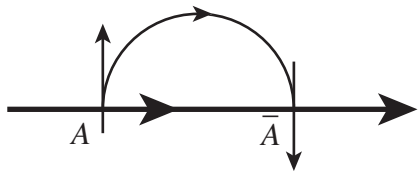}$$
    \caption{Orientations at chord endpoints} \label{F:orient}
    \end{figure}

\begin{figure} [h]
    $$\includegraphics{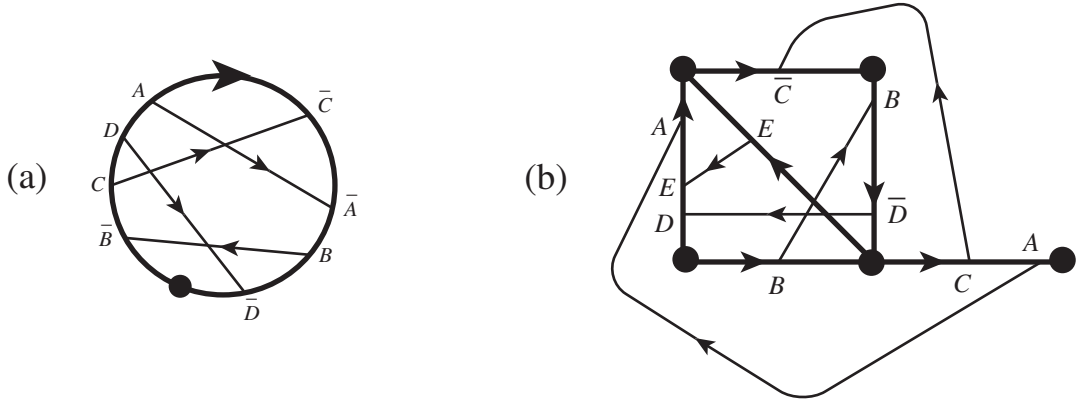}$$
    \caption{Examples of chord diagrams} \label{F:chordexamples}
    \end{figure}

Generally, we will represent a chord diagram visually by drawing an arc between the endpoints.  In an oriented chord diagram, the arc will be on the {\it same} side of the edges at the endpoints when the endpoints have {\it opposite} orientation, and on opposite sides when the endpoints have the same orientation.  This terminology is geometrically motivated - an oriented line segment in the plane naturally has left and right sides, and when the arc is on the same side of the edges at each endpoint, the bases for $\mathbb{R}^2$ given by the ordered pair of the tangent vectors to the chord and the edge have opposite orientations at the two endpoints, as shown in Figure~\ref{F:orient}.
Two examples of (oriented) chord diagrams are shown in Figure~\ref{F:chordexamples}.

\subsection{Intersection graphs for chord diagrams} \label{SS:IGchord}

Much of the work done on chord diagrams for circles has been in terms of their intersection graphs, also known as {\it circle graphs}.  These graphs have also been useful in the recent applications of chord diagrams to knot and link invariants \cite{bg, cdl, la, me1, me2}.  So it is natural to look at some kind of intersection graph for chord diagrams on more general graphs.  For chord diagrams on the circle, two chords are said to intersect if their endpoints alternate around the circle.  In a more general graph, there may not be a way to move ``around" the graph, or there may be many, so we need a different criterion.  Another point of view of circle graphs is as ``proper circular-arc graphs" - in other words, two chords intersect if the arcs on the circle between their endpoints overlap, but neither is a proper subset of the other \cite{mm}.  By replacing ``arcs on the circle" with ``paths in the graph," we can generalize this notion to other graphs.

\begin{defn} \label{D:IGchord}
Given a chord diagram C on a (possibly directed) graph G, the {\bf intersection graph for C} is the graph $\Gamma_G(C)$ whose vertices are the chords of C, and vertices $c_i$ and $c_j$ (corresponding to chords $c_i$ and $c_j$) are adjacent if given any (undirected) path $p_i$ in G between the endpoints of $c_i$ and any (undirected) path $p_j$ in G between the endpoints of $c_j$, $p_i$ and $p_j$ properly intersect (meaning that they overlap - perhaps only in one vertex - but neither is a subset of the other).
\end{defn}

\begin{figure} [h]
    $$\includegraphics{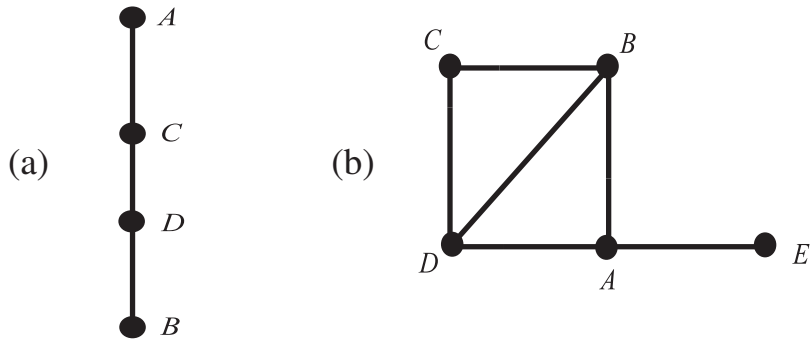}$$
    \caption{Intersection graphs for chord diagrams} \label{F:IGchord}
    \end{figure}

So for two chords to intersect, {\it all} paths between their endpoints must intersect; for them to be disjoint, there need only be one pair of paths which are disjoint (or in which one properly contains the other).  This is clearly a fairly restrictive notion of intersection, but it seems to be the one which best matches our visual intuition.  Figure~\ref{F:IGchord} shows the intersection graphs for the chord diagrams in Figure \ref{F:chordexamples}.  Note that, for chord diagrams on circles, this definition agrees with the usual definition of the intersection graph.

\section{Chord Diagrams on Planar Graphs} \label{S:planar}

A common situation is to study chord diagrams when the chords are ``real" chords - in other words, they are chords of regions in some planar embedding of the graph.  In this section, we will explore chord diagrams and their intersection graphs when we are given a planar graph together with a particular embedding in the plane.  Our main result is to give conditions when the embedding of the graph can be extended to an embedding of the entire chord diagram.

\begin{defn} \label{D:respect}
Let $G$ be a directed planar graph, $C$ be an oriented chord diagram on $G$, and $f:G \rightarrow \mathbb{R}^2$ be a particular embedding of $G$.  Consider a chord $c_i$ of $C$, with endpoints on edges $e_1$ and $e_2$ in $G$.  We say $c_i$ {\bf respects} $f$ if:
\begin{enumerate}
   \item Both endpoints of $c_i$ lie on the boundary of a single region R of $f(G)$ (i.e. a single connected component of $\mathbb{R}^2 - f(G)$).  
   \item If R is on the same side (left or right) of $e_1$ and $e_2$, then the endpoints of $c_i$ have opposite orientations; if R is on the opposite sides of $e_1$ and $e_2$, then the endpoints of $c_i$ have the same orientation (assuming $e_1$ and $e_2$ are not cut edges).
   \item If either $e_1$ or $e_2$ is a cut edge, then the endpoints of $c_i$ may have either orientation.
\end{enumerate}
If $G$ is not directed, then only the first condition is required.  We say that $C$ {\bf respects} $f$ if every chord of $C$ respects $f$.  
\end{defn}

We should observe that if $e_1$ and $e_2$ are both on the boundaries of two regions $R_1$ and $R_2$, then the definition does not depend on which region is used; if $R_1$ is on the same side of both edges, then so is $R_2$, and vice-versa.  We define the intersection graph of $C$ with respect to $f$ by restricting our attention to paths in $G$ which lie on the boundary of a single region of $f(G)$.

\begin{defn} \label{D:IGplanar}
Given a chord diagram C on a (possibly directed) planar graph G, and an embedding $f:G \rightarrow \mathbb{R}^2$, the {\bf intersection graph for C with respect to f} is the graph $\Gamma_G(C; f)$ whose vertices are the chords of C, and vertices $c_i$ and $c_j$ (corresponding to chords $c_i$ and $c_j$) are adjacent if given any (undirected) path $p_i$ on the boundary of a single region of $f(G)$ between the endpoints of $c_i$ and any (undirected) path $p_j$ on the boundary of a single region of $f(G)$ between the endpoints of $c_j$, $p_i$ and $p_j$ properly intersect.
\end{defn}

\noindent{\sc Remarks:}
\begin{enumerate}
	\item If $c_i$ does not have both endpoints on the boundary of a single region of $f(G)$, then there is no path between the endpoints which lies on the boundary of a single region.  In this case, the condition for adjacency is vacuously true, so $c_i$ is adjacent to every other chord in the chord diagram.
	\item The intersection graphs shown in Figure \ref{F:IGchord} are also the intersection graphs with respect to the embeddings shown. 
	\item  Since the number of paths in a graph may be large, it is not clear whether $\Gamma_G(C)$ can be constructed in polynomial time.  However, $\Gamma_G(C; f)$ certainly can, since an endpoint for a chord can lie on the boundary of at most two regions, there are at most 4 possible paths between the endpoints for each chord.  So the number of pairs of paths which need to be compared is at most $16\binom{n}{2}$, where $n$ is the number of chords.
	\item  If two chords intersect by Definition \ref{D:IGchord}, then they certainly intersect by Definition \ref{D:IGplanar}, but the converse is not true.  So $\Gamma_G(C) \subset \Gamma_G(C; f)$ for every embedding $f$.
\end{enumerate}

The last remark gives rise to the following open question:

\begin{quest}
Given a planar graph G and a chord diagram C on G, is $\Gamma_G(C) = \bigcap_{f: G \hookrightarrow \mathbb{R}^2}{\Gamma_G(C; f)}$?
\end{quest}

Our main focus in this section is to study when a planar embedding of a graph $G$ can be extended to give a planar embedding of a chord diagram $C$ on $G$.  To be precise, by an ``embedding of the chord diagram" we mean an embedding of the graph $H_C$ constructed from $G$ by adding vertices at the endpoints of all the chords in $C$, and edges between the endpoints of each chord in $C$.  We will also look at embeddings of {\it oriented} chord diagrams - in this case, we will insist that the embedded chords be on the same or opposite sides of the edges at each endpoint, as dictated by the orientations (see Section \ref{SS:defchord}).  Our first result is well-known (an equivalent form can be extracted from \cite{fo}, for example), but is given here for completeness and as a warm-up for our later results.

\begin{prop} \label{P:circle}
Let G be a directed graph consisting of a single loop, and let $f:G \rightarrow \mathbb{R}^2$ be an embedding of G (so $f(G)$ is a circle).  Let C be an oriented chord diagram on G.  Then f extends to an embedding of C if and only the following conditions are met:
\begin{enumerate}
	\item  $\Gamma_G(C)$ is bipartite.
	\item  The endpoints of every chord in C have opposite orientation.
\end{enumerate}
If C is not oriented, then the first condition is sufficient.
\end{prop}
{\sc Proof:}  We first note that, since every embedding of the circle in the plane is isotopic, the particular embedding $f$ is irrelevant.  And, since $f(G)$ has exactly two regions, both bounded by the circle, every path on the circle is also a path in the boundary of a region, so $\Gamma_G(C; f) = \Gamma_G(C)$.  We first show the sufficiency of the conditions of the proposition.  Since the endpoints of the chords have opposite orientations, the chord is on the same side of the circle at each endpoint, so each chord can be drawn entirely within one region of $f(G)$ (if $C$ is not oriented, we can simply choose to draw the chords on the same side of the circle at each endpoint).  Since $\Gamma_G(C)$ is bipartite, we can divide the chords into two sets, denoted the ``blue" chords and the ``red" chords, so that no two blue chords intersect, and no two red chords intersect.  We will draw the blue chords in one region of $f(G)$ (inside the circle) and the red chords in the other region (outside the circle).  Choose one of the blue chords, $c$, and connect the endpoints by an arc inside the circle.  This gives an embedding of $G \cup c$, which divides the inside of the circle into two regions.  Since no other blue chord intersects $c$, any other chord $d$ must have both endpoints on the boundary of one of these regions, and so there is an arc connecting these endpoints in $\mathbb{R}^2 - (G \cup c)$.  Continuing in this way, we can embed all of the blue arcs inside the circle; by the same argument, we can embed all the red arcs outside the circle.  This gives us an extension of the embedding to $C$.

Conversely, if we have an embedding of $C$, then each chord must lie in a single region of $f(G)$, and so must be on the same side of the circle at each endpoint, so the endpoints must have opposite orientations.  In addition, the chords inside the circle must all be disjoint, as are the chords outside the circle.  So the vertices of the intersection graph can be divided into two independent sets, and so $\Gamma_G(C)$ is bipartite.  $\Box$\\
\\
\noindent{\sc Remarks:}  
\begin{enumerate}
   \item Since there are polynomial-time algorithms to determine whether a graph is bipartite, Proposition \ref{P:circle} shows that we can determine whether the embedding extends to the chord diagram in polynomial time.
   \item Proposition \ref{P:circle} easily extends to general cycles (in which the edges may not all be oriented the same way) by replacing the condition that the endpoints of each chord have opposite orientation with the condition that every chord must respect the embedding, as in Definition \ref{D:respect}.
\end{enumerate}

\subsection{Chord diagrams on $\theta_n$-graphs} \label{SS:theta}

\begin{figure} [h]
    $$\includegraphics{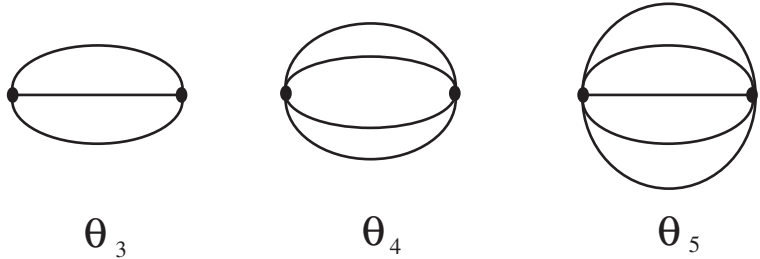}$$
    \caption{Examples of $\theta_n$-graphs} \label{F:theta}
    \end{figure}

Now we want to consider more complicated graphs.  Topologically, the simplest graphs after the circle are the $\theta_n$-graphs, which are the graphs consisting of two vertices and $n$ edges between these vertices.  Figure~\ref{F:theta} shows several examples.
The spatial embeddings of these graphs have been studied extensively \cite{go, kino, nikk}, so they are a natural set of graphs for us to consider.

\begin{prop} \label{P:theta}
Let G be a directed $\theta_n$-graph with vertices v and w, and with all edges directed from v to w.  Let $f:G \rightarrow \mathbb{R}^2$ be a planar embedding of G, and let C be an oriented chord diagram on G.  Label the edges of G by $e_1, e_2,...,e_n$ according to their (counterclockwise) cyclic order around v in f(G), and label the region of f(G) enclosed by $e_i$ and $e_{i+1}$ by $R_i$ ($R_n$ is the region enclosed by $e_n$ and $e_1$).  Then f extends to an embedding of C if and only if the following conditions are met:
\begin{enumerate}
   \item  The vertices of $\Gamma_G(C; f)$ can be labeled by the regions $R_i$ so that:
   \begin{itemize}
      \item  If a vertex is labeled $R_i$, then the corresponding chord in C has both endpoints on the boundary of region $R_i$ in f(G).
      \item  For each $R_i$, the set of vertices labeled $R_i$ is independent.
   \end{itemize}
   \item  If both endpoints of a chord are on the same edge of G, then they have opposite orientations; if the endpoints of a chord are on different edges of G, they have the same orientation (i.e. $C$ respects $f$, as in Definition \ref{D:respect}).
\end{enumerate}
If C is not oriented, then we can ignore the second condition.
\end{prop}
{\sc Proof:}  We will first show that the conditions are sufficient.  The idea is simply to draw the chord $c$ corresponding to a vertex $v_c$ labeled $R_i$ in the region $R_i$.  Since both endpoints of $c$ lie on the boundary of $R_i$, they can be connected by an arc in $R_i$.  If both endpoints of $c$ are on one edge ($e_i$ or $e_{i+1}$), then they have opposite orientations (since $C$ respects $f$), and the chord is on the same side of the edge at both endpoints; therefore, the arc is an embedding of the oriented chord.  On the other hand, if $c$ has one endpoint on $e_i$ and the other on $e_{i+1}$, then its endpoints have the same orientation, and $c$ should be on opposite sides of the edges.  But the region $R_i$ is on the left side of $e_i$ and the right side of $e_{i+1}$, so again the arc in $R_i$ is an embedding of the oriented chord.  So the embedding of $G$ can be individually extended to each chord in $C$.

We still need to show it can be extended to all the chords at once.  The obstruction to this extension is when chords cross each other - i.e. when their endpoints alternate around the boundary of a region.  We will first consider what happens if there {\it is} such an obstruction.

\begin{lem} \label{L:thetaalternate}
If two chords $c$ and $d$ both have their endpoints on the boundary of region $R_i$, and their endpoints alternate along the (undirected) cycle around the boundary of $R_i$, then the corresponding vertices $v_c$ and $v_d$ are adjacent in $\Gamma_G(C; f)$.  
\end{lem}
{\sc Proof of Lemma:}  If $c$ has both endpoints on the same edge, say $e_i$, then $d$ has exactly one endpoint on the segment $e_i$ between the endpoints of $c$, and any path between the endpoints of $d$ must contain one of the endpoints of $c$, but not both.  So any path between the endpoints of $d$ will intersect, but not contain, any path between the endpoints of $c$.  On the other hand, if neither $c$ nor $d$ have both endpoints on the same edge, then they each have one endpoint on $e_i$ and one on $e_{i+1}$.  Then any path between the endpoints of the chords with lies on the boundary of a region of $f(G)$ must run along the boundary of $R_i$, since it is the only region whose boundary contains both endpoints.  There are two such paths between the endpoints of $c$, each of which contains one endpoint of $d$.  So we conclude that if the endpoints of $c$ and $d$ alternate around the boundary of $R_i$, then $v_c$ and $v_d$ are adjacent.  This proves our lemma. $\Box$

Equivalently, if $v_c$ and $v_d$ are {\it not} adjacent, then the endpoints of $c$ and $d$ must {\it not} alternate around the boundary of $R_i$.  Since the set of vertices labeled $R_i$ is independent, none of the corresponding chords have alternating endpoints around the boundary of $R_i$.  So these chords can be embedded in $R_i$ as disjoint arcs, as in Proposition \ref{P:circle}.  This can be done for each label $R_i$, resulting in an embedding of $C$ which extends $f$.

It remains to show that the conditions are necessary.  Assume that $f$ can be extended to $C$.  Then each chord is embedded as an arc in some region $R_i$.  If both endpoints of the chord are on the same edge, the chord is on the same side at both endpoints, and the endpoints must have opposite orientations.  On the other hand, if they are on distinct edges $e_i$ and $e_{i+1}$, then the chord is on opposite sides of the edges at the two endpoints (the left side of $e_i$ and the right side of $e_{i+1}$, so the endpoints must have the same orientation.  So the second condition is necessary.

Now label the vertex $v_c$ in $\Gamma_G(C; f)$ corresponding to chord $c$ by the region in which $c$ is embedded in $f(C)$.  Obviously, both endpoints of $c$ must be on the boundary of this region.  Since we have embedded our chord diagram, chords which are in the same region $R_i$ do not cross, so their endpoints do not alternate around the boundary of $R_i$.  Then there are clearly non-overlapping paths between the endpoints of the chords, so the corresponding vertices are not adjacent in $\Gamma_G(C; f)$.  So vertices with the same label will not be adjacent.  This shows that the first condition is necessary, and completes the proof.  $\Box$\\
\\
\noindent{\sc Remark:}  It is easy to check that these conditions can be checked in polynomial time.  The coloring of most chords is forced by the location of their endpoints.  Those that are left have both endpoints on the same edge, and it is enough to show that the subgraph in $\Gamma_G(C; f)$ induced by the remaining chords on each edge is bipartite.

\subsection{Graphs without cut edges} \label{SS:nocut}

The arguments of Section \ref{SS:theta} can easily be generalized to any planar graph without {\it cut edges} (so that each edge in a planar embedding of the graph bounds a different region on either side).  Recall that $e$ is a {\it cut edge} of a connected graph $G$ if $G-e$ is not connected.

\begin{prop} \label{P:nocut}
Let G be a directed planar graph with no cut edges.  Let $f:G \rightarrow \mathbb{R}^2$ be a planar embedding of G, and let C be an oriented chord diagram on G.  Label the regions of f(G) by $R_i$ ($1\leq i \leq$ number of regions).  Then f extends to an embedding of C if and only if the following conditions are met:
\begin{enumerate}
   \item  C respects f, as in Definition \ref{D:respect} (so every chord in C has both endpoints on the boundary of some region $R_i$).
   \item  The vertices of $\Gamma_G(C; f)$ can be labeled by the regions $R_i$ so that:
   \begin{itemize}
      \item  If a vertex is labeled $R_i$, then the corresponding chord in C has both endpoints on the boundary of region $R_i$ in f(G).
      \item  For each $R_i$, the set of vertices labeled $R_i$ is independent.
   \end{itemize}
\end{enumerate}
\end{prop}

Before we begin the proof of Proposition \ref{P:nocut}, we need a technical lemma, which generalizes Lemma \ref{L:thetaalternate} in the last section.

\begin{lem} \label{L:alternate}
Consider G, f and C as in the statement of Proposition \ref{P:nocut}.  If c and d are chords in C whose endpoints alternate around the boundary of region $R$ in f(G), then the corresponding vertices $v_c$ and $v_d$ are adjacent in $\Gamma_G(C; f)$.
\end{lem}
{\sc Proof:}  Since $G$ does not have cut edges, the boundary of $R$ is a {\it cycle} in $G$.  Assume that $v_c$ and $v_d$ are {\it not} adjacent.  Then there are paths $p_c$ between the endpoints of $c$ and $p_d$ between the endpoints of $d$, each along the boundary of a region in $f(G)$, such that either $p_c \cap p_d = \emptyset$ or one path is contained in the other - without loss of generality, $p_c \subset p_d$.  So neither of the endpoints of chord $d$ lie on $p_c$.  Since the endpoints of $c$ and $d$ alternate around the boundary of $R$, $p_c$ cannot be contained in the boundary of $R$, or it would contain one of the endpoints of $d$.  So $p_c$ is in the boundary of some other region $S$, meaning that the endpoints of $c$ lie on the boundary between regions $R$ and $S$.

\begin{figure} [h]
    $$\includegraphics{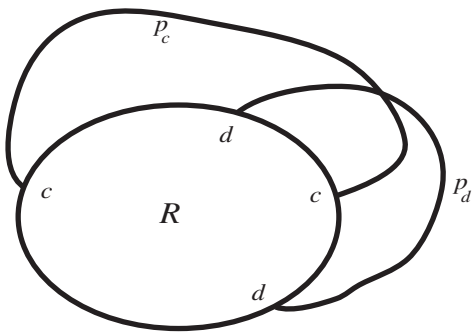}$$
    \caption{Figure for Lemma \ref{L:alternate}} \label{F:alternate}
    \end{figure}

First consider the case when $p_c$ and $p_d$ are disjoint.  Then the union of $p_c$ with one of the paths between the endpoints of $c$ along the boundary of region $R$ gives a cycle $\gamma$ (or at least a circuit which contains a cycle $\gamma$) which contains one endpoint of $d$.  So $p_d$ must intersect $\gamma$.  If $p_d$ lies along region $R$, then it will contain one of the endpoints of $c$, and so won't be disjoint from $p_c$.  But since $p_d$ can't pass {\it through} region $R$ (since this is a region of the planar embedding of $G$), it must (by the Jordan Curve Theorem) pass from one side of $\gamma$ to the other, and so must intersect $p_c$, as shown in Figure~\ref{F:alternate}.
This contradicts the assumption that $p_c$ and $p_d$ are disjoint.

Therefore, $p_c \subset p_d$, so the endpoints of $c$ are in $p_d$.  Since the endpoints of $c$ are on the boundary between regions $R$ and $S$, and we know that $p_d$ cannot be on the boundary of region $R$, it must lie on the boundary of region $S$.  So the endpoints of $d$ are also on the boundary between regions $R$ and $S$.  Since $p_c \subset p_d$, the endpoints of $c$ and $d$ do {\it not} alternate around the boundary of region $S$, which is a cycle.  So there are paths $p_c'$ and $p_d'$ between the endpoints of chords $c$ and $d$ on the boundary of region $S$, with $p_c' \cap p_d' = \emptyset$.  But this is impossible, by the argument in the last paragraph.

Therefore, $p_c$ and $p_d$ must properly intersect, so $v_c$ and $v_d$ are adjacent. $\Box$\\
\\
\noindent{\sc Proof of Proposition \ref{P:nocut}:}  Our proof is very similar to the proof of Proposition \ref{P:theta}.  We will first show that the conditions are sufficient.  Once again, the idea is simply to draw the chord $c$ corresponding to a vertex $v_c$ labeled $R$ in the region $R$.  Since both endpoints of $c$ lie on the boundary of $R$, they can be connected by an arc in $R$.  The first condition ensures that this arc will be an embedding of the oriented chord $c$.  So the embedding of $G$ can be individually extended to each chord in $C$.

We still need to show it can be extended to all the chords at once.  Let $d$ be another chord on region $R$, with corresponding vertex $v_d$.  By Lemma \ref{L:alternate}, if $v_c$ and $v_d$ are {\it not} adjacent, then the endpoints of $c$ and $d$ must {\it not} alternate around the boundary of $R$.  Since the set of vertices labeled $R$ is independent, none of the corresponding chords have alternating endpoints around the boundary of $R$.  So these chords can be embedded in $R$ as disjoint arcs, as in Proposition \ref{P:circle}.  This can be done for each label $R$, resulting in an embedding of $C$ which extends $f$.

It remains to show that the conditions are necessary.  Assume that $f$ can be extended to $C$.  Then each chord is embedded as an arc in some region $R$, with endpoints on edges $e_1$ and $e_2$.  If $R$ is on the same side of both edges, then so is the embedded chord, and the endpoints must have opposite orientations; if it is on opposite sides, the endpoints must have the same orientation.  So the first condition is necessary.

Now label the vertex $v_c$ in $\Gamma_G(C; f)$ corresponding to chord $c$ by the region in which $c$ is embedded in $f(C)$.  Obviously, both endpoints of $c$ must be on the boundary of this region.  Chords which are in the same region $R$ do not cross, so their endpoints do not alternate around the boundary of $R$ (since the boundary of $R$ is a cycle).  Then there are clearly non-overlapping paths between the endpoints of the chords, so the corresponding vertices are not adjacent in $\Gamma_G(C; f)$.  So vertices with the same label will not be adjacent.  This shows that the second condition is necessary, and completes the proof.  $\Box$

\begin{cor} \label{C:nocut}
Given a planar graph $G$ with no cut edges, a chord diagram $C$ on $G$, and a planar embedding $f:G \rightarrow \mathbb{R}^2$, there is a polynomial time algorithm to determine whether $f$ extends to an embedding of $C$.
\end{cor}

\noindent{\sc Proof:}  The conditions of Proposition \ref{P:nocut} can be checked in polynomial time.  The coloring of most chords is forced by the location of their endpoints.  Those that are left have both endpoints on boundaries of the same two regions, and it is enough to show that the subgraph in $\Gamma_G(C; f)$ induced by the remaining chords lying on each pair of regions is bipartite.  $\Box$

\subsection{Chord diagrams on line segments} \label{SS:chordline}

\begin{figure} [h]
    $$\includegraphics{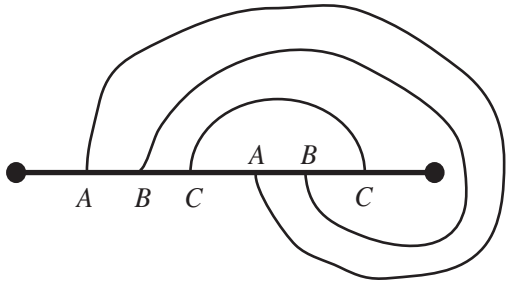}$$
    \caption{A chord diagram on $K_2$} \label{F:lineexample}
    \end{figure}

Now we want to consider planar graphs which {\it do} have cut edges.  As an introduction, we will consider the simplest such graph, the complete graph on two vertices $K_2$ (i.e. the line segment consisting of two endpoints and the edge between them).  Why do cut edges pose difficulties?  One difficulty is that the boundaries of the regions in a planar embedding of the graph are no longer cycles in the graph.  For example, we might expect a line segment to be no more complicated than a circle.  But consider the example in Figure~\ref{F:lineexample}.
The sequence of endpoints along the graph is $ABCABC$.  According to Definition \ref{D:IGplanar}, all three of these chords intersect.  So the intersection graph has a cycle of length 3, and so is not bipartite.  And yet we {\it can} extend the embedding of the line segment to an embedding of the chord diagram, so Proposition \ref{P:circle} does not extend to $K_2$.  The fact that chords can go ``around the end" of the graph increases the difficulty of the problem.

However, we would still like to exploit our results from Section \ref{SS:nocut}.  So our plan is to transform a graph {\it with} cut edges into a graph {\it without} cut edges, and hope solving the problem in this new graph will give us the solution for the original graph.  To this end, we introduce an operation called ``blowing up" an edge in a planar embedding of a graph or chord diagram.

\begin{defn} \label{D:blowup}
Let G be a planar directed graph with an edge e (which is not a loop), and let $f:G \rightarrow \mathbb{R}^2$ be a planar embedding of G.  Say that the endpoints of e are vertices v and w, with e oriented from v to w.  Then the {\bf blowup of $e$ in $G$ with respect to $f$} is the graph $G_e$ formed by removing edge $e$ and replacing it with two edges e and $e'$, both directed from $v$ to $w$.  Figure~\ref{F:blowup} illustrates this operation.
    \begin{figure} [h]
    $$\includegraphics{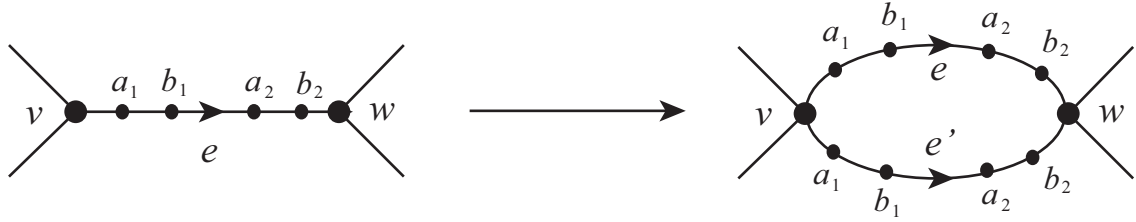}$$
    \caption{Blowing up edge $e$} \label{F:blowup}
    \end{figure}
\end{defn}
{\sc Remark:}  Figure \ref{F:blowup} also shows that the embedding $f$ of $G$ is easily modified to give an embedding $f_e$ of $G_e$.  In particular, this means that $G_e$ is also a planar graph.

What is the effect of blowing up an edge $e$ on a chord diagram $C$ on $G$?  If $a$ is a chord of $C$, label its endpoints by $a_1$ and $a_2$.  If $a_i$ is on edge $e$ in $G$, it will be doubled in $G_e$, appearing on both $e$ and $e'$.  In other words, the sequence of endpoints along $e$ and $e'$ will both be the same as the original sequence along $e$.  The ``chord" $a$ now consists of all pairings of a point labeled $a_1$ and a point labeled $a_2$; there are $2^k$ such pairings, where $k = 0, 1, 2$ is the number of endpoints of the original chord $a$ which were on edge $e$.  We give each pairing its own name, $a, a', a'', a'''$, yielding a new chord diagram $C_e$ on $G_e$.

\begin{figure} [h]
    $$\includegraphics{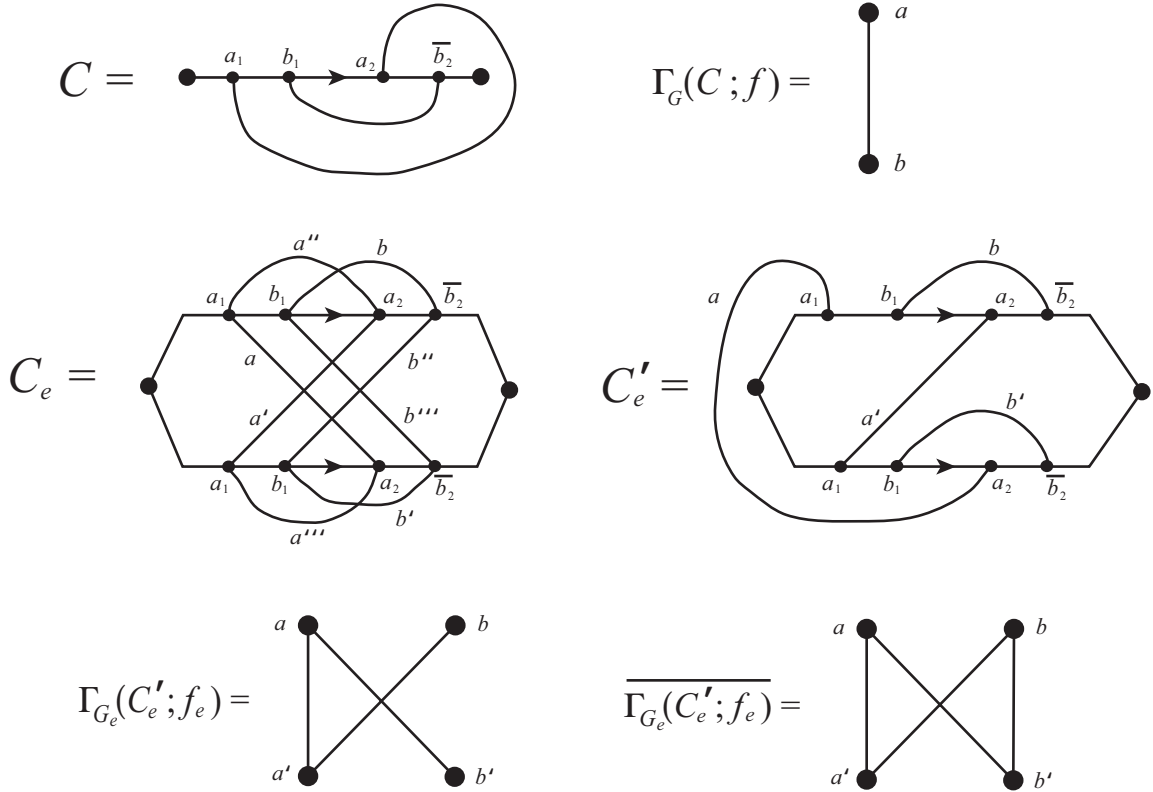}$$
    \caption{New chord diagrams obtained from blowing up an edge} \label{F:chordblowup}
    \end{figure}

If the chord diagram $C$ is oriented, then the orientations of $a_1$ and $a_2$ are also replicated in $G_e$.  If the chord $a$ in $C$ respected the embedding $f$, and $a$ had an endpoint on $e$, then only half of the set of chords labeled $a$ in $G_e$ will respect the embedding $f_e$.  So we obtain the {\it reduced} oriented chord diagram $C_e'$ on $G_e$ by discarding those chords $a$ which do not respect $f_e$.  The intersection graph $\Gamma_{G_e}(C_e', f_e)$ is defined as before.  We will add to this graph edges between all the chords which came from the chord $a$, to get a new intersection graph $\overline{\Gamma_{G_e}(C_e'; f_e)}$.  Figure~\ref{F:chordblowup} illustrates this process for an example with $G = K_2$.

\begin{prop} \label{P:line}
Let $G=K_2$ be a directed line segment, and let $f:G \rightarrow \mathbb{R}^2$ be an embedding.  Let $C$ be an oriented chord diagram on $G$.  Let $G_e$ be the result of blowing up the single edge $e$ of $G$, $f_e$ be the resulting embedding of $G_e$, and $C_e'$ be the reduced oriented chord diagram on $G_e$ obtained from $C$.  Then $f$ can be extended to an embedding of $C$ if and only if $\overline{\Gamma_{G_e}(C_e'; f_e)}$ is bipartite.
\end{prop}
{\sc Proof:}  In this case, every chord in $C$ has both endpoints on $e$, so each chord of $C$ corresponds to 4 chords in $C_e$, and hence 2 chords in $C_e'$.  We first show the sufficiency.  If  $\overline{\Gamma_{G_e}(C_e'; f_e)}$ is bipartite, then so is its subgraph $\Gamma_{G_e}(C_e', f_e)$.  Since every chord in $C$ respects $f$ (since the only edge is a cut edge), $C$ respects $f$; so by definition $C_e'$ respects $f_e$ and contains a chord corresponding to every chord in $C$.  So by Proposition \ref{P:circle}, $f_e$ extends to an embedding of $C_e'$.  Moreover, by using the partition of vertices from $\overline{\Gamma_{G_e}(C_e'; f_e)}$, we can extend $f_e$ so that exactly one of each pair of chords $a, a'$ is inside the cycle $G_e$ and one is outside.  The chords which are outside give an extension of $f$ to $C$.

Conversely, if $f$ extends to $C$, then we can draw the chords of $C$ as disjoint arcs in the single region of $f(G)$.  These arcs can then be drawn as arcs outside the cycle $G_e$ in $f_e(G_e)$.  But there is another extension of $f$ to $C$ obtained by taking the reflection of the first embedding across the axis defined by the image of $G$.  These arcs can be drawn {\it inside} the circle $f_e(G_e)$.  Together, these give an extension of $f_e$ to an embedding of $C_e'$, so $\Gamma_{G_e}(C_e', f_e)$ is bipartite.  Since the arcs for the chords $a$ and $a'$ in this extension are on opposite sides of the circle, adding edges between them in the intersection graph will not affect the partition into two independent sets, so $\overline{\Gamma_{G_e}(C_e'; f_e)}$ is also bipartite. $\Box$

Notice that, in Figure \ref{F:chordblowup}, the intersection graph $\overline{\Gamma_{G_e}(C_e'; f_e)}$ is bipartite, and the chord diagram $C_e'$ can be embedded as shown.

\subsection{Graphs with cut edges} \label{SS:cut}

We can extend the ideas of the last section to prove the analogue of Proposition \ref{P:nocut} for chord diagrams on graphs with cut edges.  Our plan is simply to blow up every cut edge, as in Section \ref{SS:chordline}, and consider the resulting graph, which no longer has cut edges.

\begin{figure} [h]
    $$\includegraphics{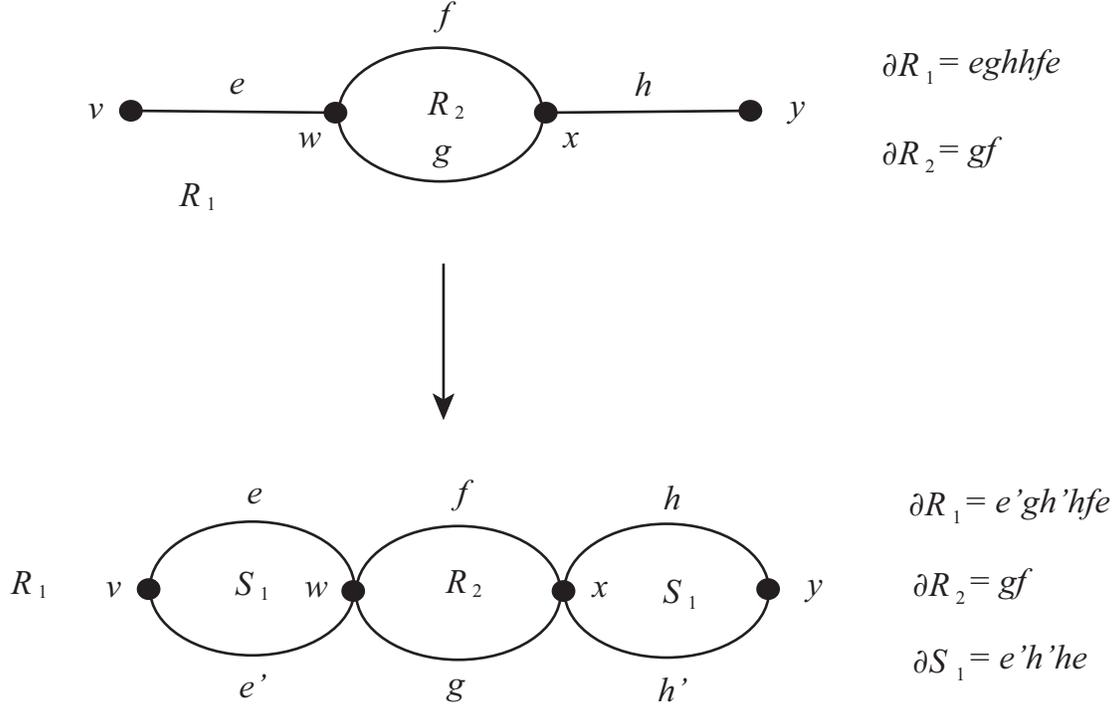}$$
    \caption{Regions $R_i$ and $S_i$ and their boundaries} \label{F:boundaries}
    \end{figure}

Consider a graph $G$, and let $E$ be the set of cut edges.  Let $G_E$ be the result of blowing up every edge in $E$.  If $C$ is an oriented chord diagram on $G$, then $C_E$ is the corresponding chord diagram on $G_E$.  If $f:G \rightarrow \mathbb{R}^2$ is an embedding of $G$ in the plane, then we have a corresponding embedding $f_E$ of $G_E$.  The regions of $f_E(G_E)$ are essentially the same as the regions of $f(G)$, together with the new regions created by the blown up edges (one new region for each cut edge in $G$).  If we label the regions of $f(G)$ by $R_1,..., R_n$, then we can give the corresponding regions of $f_E(G_E)$ the same labels; in addition, we will give {\it all} the regions created by blowing up cut edges on the boundary of $R_i$ the same label $S_i$.  The ``region" $S_i$ is the union of all these regions, and its boundary is the union of the images $e, e'$ of the cut edges of $R_i$ in $C$, in the same order that they are found in the boundary of $R_i$ in $C_E$ (so the boundary of $S_i$ is a subset of the boundary of $R_i$).  This is illustrated in Figure \ref{F:boundaries}.

As usual, a chord of $C_E$ {\it respects} $f_E$ if it satisfies Definition \ref{D:respect}, treating $S_i$ as a single region.  So we can consider the reduced chord diagram $C_E'$ which results from removing all the chords in $C_E$ which do not respect $f_E$.  

The intersection graph $\Gamma_{G_E}(C_E'; f_E)$ is defined as usual, with the exception that a ``path" on the boundary of a region $S_i$ is taken to be a sequence of consecutive edges in the boundary of $S_i$ (as defined in the last paragraph), even though these edges may not be adjacent in $G_E$.  As in Section \ref{SS:chordline}, $\overline{\Gamma_{G_E}(C_E'; f_E)}$ is defined by adding to $\Gamma_{G_E}(C_E'; f_E)$ edges between vertices corresponding to chords in $C_E'$ which came from the same chord in $C$.

\begin{prop} \label{P:cut}
Let G be a directed planar graph with a set of cut edges E, and let $G_E$ be the result of blowing up the edges in E, as described above.  Let $f:G \rightarrow \mathbb{R}^2$ be a planar embedding of G, and let C be an oriented chord diagram on G.  Label the regions of f(G) by $R_i$ ($1\leq i \leq$ number of regions).  We define $f_E$ and $C_E$ as described above.  Then f extends to an embedding of C if and only if the following conditions are met:
\begin{enumerate}
   \item  C respects f.
   \item  The vertices of $\overline{\Gamma_{G_E}(C_E'; f_E)}$ can be labeled by the regions $R_i$ and $S_i$ so that:
   \begin{itemize}
      \item  If a vertex is labeled $R$, then the corresponding chord in $C_E'$ has both endpoints on the boundary of region $R$ in $f_E(G_E)$.
      \item  For each region $R$ ($R_i$ or $S_i$), the set of vertices labeled $R$ is independent.
      \item  If $a$ was a chord of $C$ with both endpoints on a cut edge in the boundary of $R_i$, then exactly one of the chords $a, a'$ in $C_E'$ is labeled $R_i$, and the other is labeled $S_i$.
   \end{itemize}
\end{enumerate}
\end{prop}
{\sc Proof:}  We first show the sufficiency.  From Proposition \ref{P:nocut}, the conditions imply that $f_E$ extends to an embedding of $C_E'$.  Moreover, for every chord in $C$, there is exactly one corresponding chord of $C_E'$ whose image lies in the union of the regions $R_i$ (as opposed to the regions $S_i$).  So we can then ``collapse" the regions $S_i$ (reversing the blowup operations) to obtain an embedding of $C$ which extends $f$.

For the necessity, assume that $f$ does extend to $C$.  Then, for any chord $a$ of $C$ such that $f(a)$ is an arc in the region $R_i$, there is a corresponding arc in the region $R_i$ in the complement of $f_E(G_E)$.  These arcs are all still disjoint.  This provides an embedding of almost all chords of $C_E'$; label the vertices in $\overline{\Gamma_{G_E}(C_E'; f_E)}$ corresponding these chords by the region in which each is embedded.  The vertices labeled $R_i$ will still be independent, as they were in $\Gamma_G(C; f)$.  The only remaining chords are ``doubles" of chords which have already been embedded in some $R_i$, and which correspond to chords in $C$ which had both endpoints on cut edges on the boundary of $R_i$.  So these chords have both endpoints on the boundary of $S_i$, and the corresponding vertices can be labeled $S_i$.  

It only remains to check that the vertices labeled $S_i$ will be independent.  However, the order of the endpoints of the chords labeled $S_i$ around the ``boundary" of $S_i$ is exactly the reverse of the order of the endpoints of their ``doubles" around the boundary of $R_i$.  So the endpoints of two chords alternate around the boundary of $S_i$ only if the endpoints of their doubles alternate around the boundary of $R_i$.  Since the vertices labeled $R_i$ are all independent, so are the vertices labeled $S_i$.  So the conditions of the proposition are satisfied, completing the proof. $\Box$

\begin{cor} \label{C:cut}
Given a planar graph $G$, a chord diagram $C$ on $G$, and a planar embedding $f:G \rightarrow \mathbb{R}^2$, there is a polynomial time algorithm to determine whether $f$ extends to an embedding of $C$.
\end{cor}

\noindent{\sc Proof:}  From Corollary \ref{C:nocut}, we can check whether $f_E$ extends to $C_E'$ in polynomial time.  The only condition that is left to check is the last one, that one of each pair $a, a'$ can be labeled $R_i$ and the other $S_i$.  Since, for a connected graph, there is at most one division into two independent sets (vertices with paths which are odd or even length to a given vertex), this can be done very quickly after determining that the set of chords with endpoints on $R_i$ and $S_i$ is bipartite.  $\Box$

\section{Gauss Codes for Graphs} \label{S:gauss}

As an application, we will use chord diagrams and intersection graphs for graphs to study Gauss codes for graphs.  The Gauss code was first developed by C. F. Gauss to study closed curves in the plane which intersect themselves only in transverse double points (called ``crossings").  Label the crossings by some set of symbols (such as integers or the letters of the alphabet), and give the curve an orientation.  Beginning at an arbitrary crossing, we obtain a Gauss code for the curve by writing down the sequence of labels of the crossings passed as we traverse the curve, following its orientation.  The Gauss code provides a convenient combinatorial representation of the curve.  However, not all sequences which ``look like" Gauss codes actually represent closed curves in the plane.  The problem is to determine (in polynomial time) which sequences of symbols can be obtained from a closed curve in this way - i.e., which sequences are {\it realizable} as the Gauss code for a closed curve in the plane.  There are several different solutions to this problem \cite{fo,rr}, we will describe one in Section~\ref{SS:circle} which uses Proposition \ref{P:circle}.

Kauffman's recent development of the theory of virtual knots \cite{ka} was motivated in part by a desire to realize the ``unrealizable" Gauss codes.  This paper was largely motivated by a desire to extend Kauffman's work to spatial graphs, and study ``virtual" spatial graphs.  As a preliminary, we need to extend the notion of a Gauss code to general graphs, and consider when an abstract code is realizable by an immersion of a graph in a plane.  In future work, this will lead to the introduction of virtual spatial graphs \cite{fm}.

Our goal is to generalize to the situation where we have an immersion of an arbitrary graph in the plane, rather than a circle.  In this case, as described in Section~\ref{SS:graph}, the ``Gauss code" is a set of sequences associated with the edges of the graph.  We will give a algorithm, based on Proposition \ref{P:cut}, for determining whether an abstract Gauss code can be realized as the Gauss code for an immersed graph and, if so, constructing the desired immersion.

\subsection{Gauss Codes for Closed Curves} \label{SS:circle}

We will first describe a method for determining whether a classical Gauss code is realizable by an immersed
closed curve in the plane.  This problem was first solved by Dehn \cite{de}, and a clear presentation of his solution is
given by Read and Rosenstiehl \cite{rr}.  Our method builds on Dehn's, and is equivalent to that given by de Fraysseix and Ossona de Mendez \cite{fo}, but is presented in the language of chord diagrams developed in Section \ref{S:chord}.

\begin{figure} [h]
    $$\includegraphics{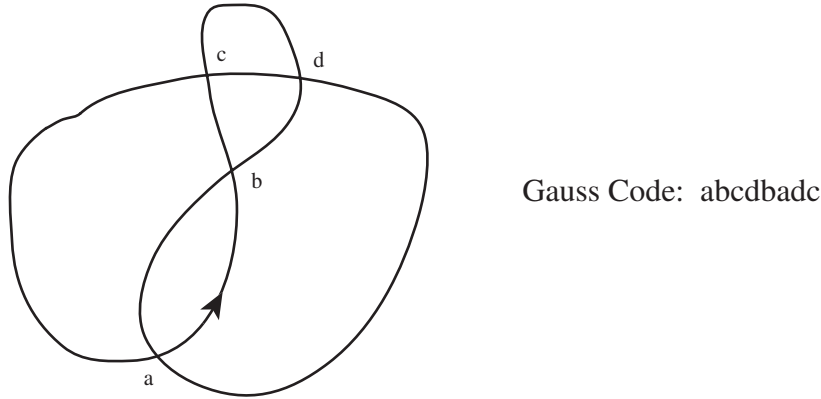}$$
    \caption{Example of a Gauss code} \label{F:gaussloop}
    \end{figure}

\begin{figure} [h]
    $$\includegraphics{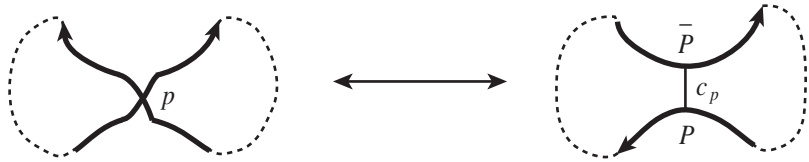}$$
    \caption{Splitting a crossing} \label{F:split}
    \end{figure}

An oriented closed curve in the plane gives rise to a crossing sequence, the {\it Gauss code} for the curve,
by labeling the $n$ self-crossings of the curve, and then writing down the sequence of $2n$ labels
encountered as we trace the curve following its orientation (obviously, the code is only unique modulo
cyclic permutations).  An example is shown in Figure~\ref{F:gaussloop}.
The problem, initially posed by Gauss, is to determine whether an arbitrary sequence of length $2n$,
containing two occurrences each of $n$ symbols, is realizable as the Gauss code of a closed plane curve. 
We will call such a sequence a {\it crossing sequence}, and refer to the symbols as {\it crossings}.  We
follow Read and Rosenstiehl by defining a {\it splitting} of a crossing sequence at each crossing, and then studying the
resulting {\it split sequence} or {\it split code}.  The notion of splitting a crossing is motivated by the
idea of smoothing a self-intersection of a curve in the plane into two non-intersecting arcs.  We consider a small neighborhood of a self-crossing $p$.  Inside this neighborhood, we replace the two intersecting arcs with two disjoint arcs with the same endpoints, as shown in Figure~\ref{F:split}.  
This replacement is done so as to preserve the number of components of the curve, which requires reversing the orientation of one of the arcs from $p$ to itself, so that the resulting curve is still oriented.  Aside from the change of orientation, the operation is entirely local, so if we label points $P$ and $\bar{P}$ on the two new arcs, we can draw a chord $c_p$ between them which does not intersect the rest of the curve and is contained within the original neighborhood of the crossing, as shown in Figure \ref{F:split}.  Notice that, if the chord $c_p$ is oriented, its endpoints at $P$ and $\bar{P}$ will have opposite orientations.

We can translate this splitting operation into moves which can be performed directly on a crossing sequence, regardless of whether or not the sequence is realizable.  

\begin{defn} \label{D:splitting}
Given a sequence $S = \alpha p \beta p \gamma$, {\bf splitting S at p} means to replace S by
$\alpha P \beta^{-1} \bar{P} \gamma$, where we say that $P$ and $\bar{P}$ have opposite orientations,
and $\beta^{-1}$ is the result of writing the subsequence $\beta$ in reverse order and reversing the
orientation of any previously split symbols in $\beta$.
\end{defn}

\begin{figure} [h]
    $$\includegraphics{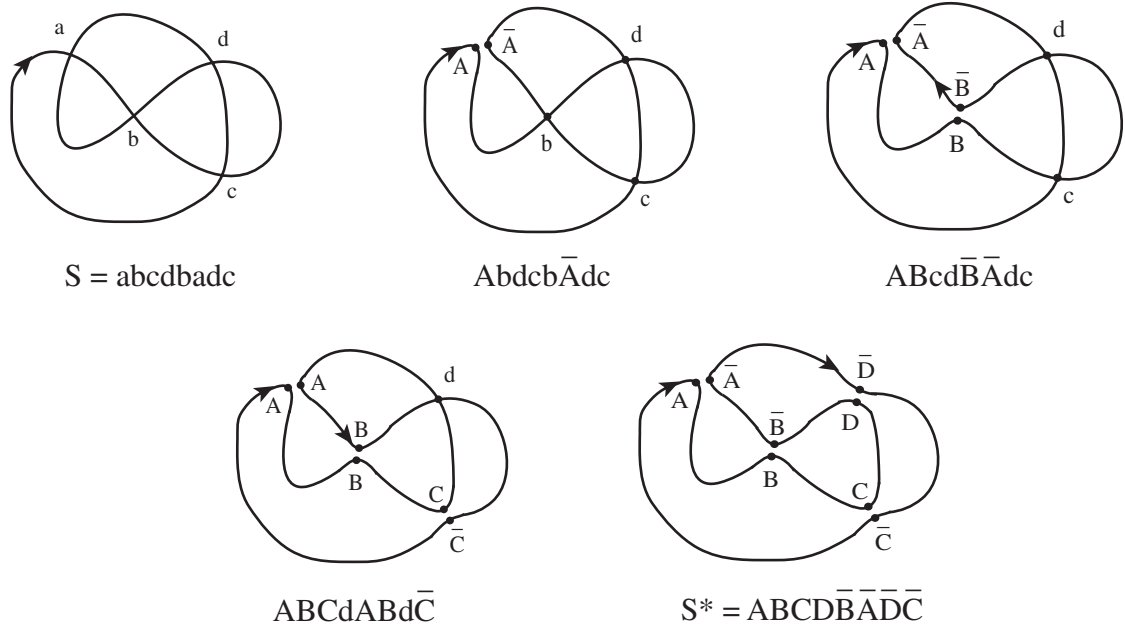}$$
    \caption{Obtaining the split code} \label{F:splitcode}
    \end{figure}

Given a crossing sequence $S$, the result of splitting $S$ at {\it every} crossing is called a {\it split
sequence} or {\it split code} for $S$, and is denoted $S^*$.  The split code is {\it not} unique - it
depends upon the order in which the crossings are split.  Figure~\ref{F:splitcode} illustrates this
procedure with an example of smoothing all the crossings in a closed plane curve, and the corresponding
operations on the Gauss code in order to obtain the split code.
The result of performing the splitting process on a closed plane curve is a {\it simple} closed curve (i.e. an embedded circle) labeled by the split code.  We will call this the {\it split curve}.  Clearly, the split curve is an oriented chord diagram on the circle, in the sense of Definition \ref{D:orientchord}.  To reconstruct the original curve from the split curve, we simply need to connect each pair of points $P$ and $\bar{P}$ by the chord $c_p$ shown in Figure \ref{F:split}, and contract the chords to bring the points back together, reversing the splitting process.  So reconstructing the original curve is done by extending the embedding of the underlying circle of the split curve to an embedding of the oriented chord diagram represented by the split curve.  This leads directly to the following result.

\begin{thm} \label{T:closedcurve}
\cite{fo, rr} A crossing sequence S is realizable as the Gauss code of a closed curve in the plane if and only if the the oriented chord diagram on the circle represented by the split sequence S* can be embedded in the plane.
\end{thm}
{\sc Proof:}  From the discussion above, if the crossing sequence $S$ is realizable, then there is an embedding of the chord diagram $S^*$.  Conversely, an embedding of $S^*$ can be turned into a collection of intersecting closed loops by expanding each chord into a crossing, reversing the splitting operation in Figure \ref{F:split}.  It only remains to check that this process will yield a single closed curve with Gauss code $S$ (this is the {\it complete traceability} of \cite{rr}).

Say that $S^*$ was produced from $S$ via a series of sequences $S = S_0, S_1, S_2, ..., S_n, S_{n+1} = S^*$, where $S_{i+1}$ is the result of splitting crossing $p_i$ in $S_i$.  We will expand the chords between pairs of symbols in $S^*$ in the reverse order from how they were split - i.e. the last symbol to be split is the first chord added.  So we begin with the chord $c_{p_n}$ between $P_n$ and $\bar{P_n}$ in $C_{n+1} = C^*$.  Expanding this chord into a crossing, and reversing the orientation on the arc of $C_{n+1}$ originally directed from $P_n$ to $\bar{P_n}$, exactly reverses the splitting process of Definition \ref{D:splitting}.  A key observation is that, since $P_n$ and $\bar{P_n}$ have opposite orientations, the result of expanding the chord is a closed curve with one component.  Since we have exactly reversed the final splitting operation, the result is an oriented closed curve $C_n$ whose Gauss code is $S_n$.  We can continue this process, always choosing to expand the chord corresponding to the last crossing which was split, so that its endpoints have opposite orientation and the resulting curve still has a single component.  In this way, we can successively construct an oriented closed curve $C_i$ for each sequence $S_i$.  Ultimately, we will have a curve $C$ whose Gauss code is the original crossing sequence $S$, as desired.  $\Box$\\
\\
\noindent {\sc Remarks:}  There are several remarks that should be made at this point.
\begin{enumerate}
     \item To check whether a particular split code is embeddable, we apply the criteria of Proposition \ref{P:circle}.  So we check that the two occurrences of each symbol have opposite orientation and that the intersection graph is bipartite.  The importance of Theorem \ref{T:closedcurve} is that this can be done in polynomial time.
     \item Comparing with the {\it D-switch} operation of de Fraysseix and Ossona de Mendez \cite{fo}, the condition that the two occurrences of each symbol have opposite orientation is equivalent to saying that a symbol and its twin are not adjacent in the interlacement graph.  So, the result of Theorem \ref{T:closedcurve} then also follows from the proof of Theorem 6 in \cite{fo}.
     \item Kauffman \cite{ka} provides another characterization, in which the condition that each pair of symbols have opposite orientations is replaced by the condition that the crossing sequence $S$ is {\it evenly intersticed}, meaning that between two occurrences of any symbol is a word of even length.  However, Kauffman does not prove that this requirement is sufficient.  It is not hard to show that, if $S$ is evenly intersticed, then the number of unsplit symbols between a pair of split symbols (at any stage in the process of splitting $S$) is even exactly when the split symbols have opposite orientations.  Since, in $S^*$, the number of unsplit symbols is 0, and hence even, all pairs of split symbols have opposite orientations, so by Proposition \ref{P:circle} and Theorem \ref{T:closedcurve} the crossing sequence is realizable.  This shows that Kauffman's characterization is also valid.
     \item Not every sequence of symbols $S^*$ which satisfies the conditions of Theorem \ref{T:closedcurve} is the split code of a crossing sequence.  For example, if the reconstruction algorithm is applied to the sequence $AB\bar{A}\bar{B}$, the result will be a curve with two components; however, this is not a contradiction because this sequence {\it cannot} be obtained as the split code of a crossing sequence.  To see this, observe that if we ``unsplit" symbol $B$ we are left with $AbAb$ - since the two occurrences of $A$ have the same orientation, this could not be the result of a splitting operation; attempting to unsplit $A$ leads to a similar situation.
\end{enumerate}

\subsection{Gauss Codes for Graphs} \label{SS:graph}

We would like to extend the ideas of the last section to general directed graphs.  In this section, our graphs will be {\it
vertex-oriented}, meaning that we have a cyclic order on the edges incident to each vertex.

Given a graph $G$, a {\it Gauss code for G} consists of a set of symbols $A$ and a set of sequences in
these symbols, one for each edge of the graph.  Each sequence will begin with the vertex at the start of the
edge, and the position of the edge in the cyclic ordering around that vertex, and end with the vertex at
the end of the edge, with the position of the edge in the cyclic ordering around the end vertex.  In
between are symbols from $A$.  Notice that the abstract graph $G$ can be reconstructed from the set of Gauss
codes.

\begin{figure} [h]
    $$\includegraphics{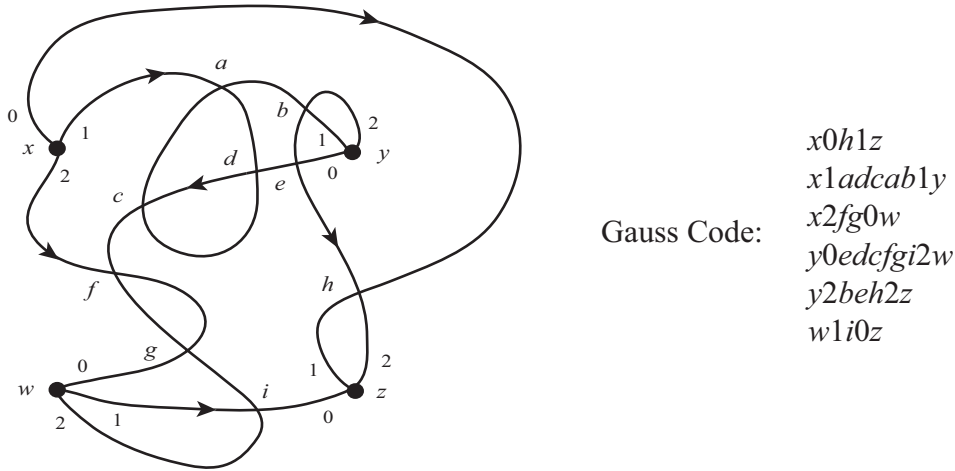}$$
    \caption{Example of a Gauss code from a graph immersion} \label{F:example}
    \end{figure}

While an abstract graph has many possible Gauss codes, we can assign a unique code (up to cyclic permutation of the vertex orderings) to any immersion of the graph in the plane, just as we do for immersions of circles.  In this case, the sequence of symbols just records the order of the crossings along each edge.  An example is shown in Figure~\ref{F:example}.
We wish to address the problem of determining which Gauss codes for $G$ can be realized by an immersion of
$G$ in the plane.  We will assume that the initial graph $G$ is {\it connected}.

\begin{figure} [h]
    $$\includegraphics{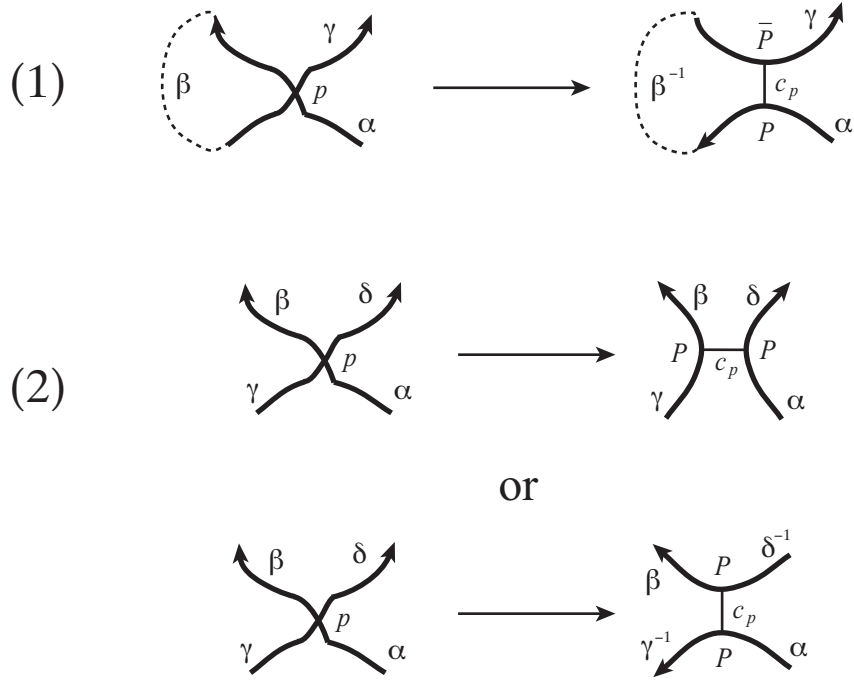}$$
    \caption{Ways of splitting a crossing} \label{F:split2}
    \end{figure}

Our first observation is that every symbol in the set of sequences (other than the vertices) must appear exactly twice.  This allows us to create a {\it split code} from the original Gauss code.  As in the previous sections, the motivation for the split code is to split the crossings of an immersion of a graph in the plane.  There are two ways to split a crossing (one is shown in Figure~\ref{F:split}) - we want to perform the split so that the graph remains connected.  If both splits preserve the connectedness of the graph, we can choose one arbitrarily.  As before, we also want to keep track of the relative local orientations of the smoothed arcs.  These operations can be defined abstractly in terms of the Gauss code.  If $p$ is a symbol in the alphabet $A$, then splitting the code at $p$ means one of the following operations:
\begin{enumerate}
   \item  If there is a sequence $w = \alpha p\beta p\gamma$, then the new sequence is $w^* = \alpha
P\beta^{-1} \bar{P}\gamma$.  So the subsequence $\beta$ has been reversed.
   \item  If there are two sequences $w = \alpha p \beta$ and $u = \gamma p \delta$, then the new sequences
are  either $w^* = \alpha P \delta$ and $u^* = \gamma P \beta$ or $w^* = \alpha P \gamma^{-1}$ and $u^* =
\delta^{-1} P \beta$.
\end{enumerate}
These are illustrated in Figure~\ref{F:split2}.  Notice that, in the second case, the underlying graph is changed in either choice, but at least one of these graphs is connected.  We will always make a choice which keeps the underlying graph connected.  Also notice that the first choice preserves orientations along each edge segment, while the second choice reverses the orientations along segments coming from one of the two edges.

The result of performing this splitting operation on every symbol in the set of sequences gives a split code for the original Gauss code.  The underlying abstract graph for the split code is called the {\it split graph}.  The splitting process is not unique - the final set of sequences after all crossings have been split will depend on the order in which they were split, and how they were split.  Any of these possibilities is called a split code for the original Gauss code.  A split code can be viewed as an oriented chord diagram on its split graph.

\begin{lem} \label{L:planar}
If a Gauss code is derived from an immersion of a connected graph, then the underlying split graph for any split
code is planar and connected; moreover, the split code determines a unique (up to isotopy) spherical
embedding of the split graph.
\end{lem}
{\sc Proof:}  If we begin with an immersion of a graph and split every crossing, we are left with a graph
embedded in the plane; this implies that the graph underlying the split code is planar, and the split code
is realized by some planar embedding.  The embedding is connected since we choose our splittings to keep
the graph connected.  It remains to show that this embedding is unique (up to isotopy on the sphere).

We will describe an algorithm for constructing a planar embedding from the split code.  We will show that
(up to isotopy on the sphere) the algorithm produces at most one embedding, and that it fails to produce an
embedding respecting the split code only if one does not exist.  In fact, we will not use all the
information of the split code, only the underlying graph and the order of the edges around each vertex. 
The idea of the algorithm is to trace out faces of the planar embedding.

Begin with an arbitrary vertex $v$, and an edge $e$ incident to $v$.  Embed $v$ and $e$.  Say that $w$ is
the other endpoint of $e$.  Let $f$ be the next edge adjacent to $w$ after $v$ (using the cyclic ordering
of the edges adjacent to $w$), and add $f$ to the embedded graph.  Continue this process, moving one edge
around each vertex.

The first time we return to a vertex we have visited before, we will have two choices of how to draw the
edge - clockwise or counterclockwise around the previously embedded parts of the graph.  But these choices
are isotopic on the sphere.  Once this first choice is made, all future edges will connect two points on
the boundary of a (topological) disk, so all possible ways to draw them will be isotopic.  As a result,
there are no more choices to be made, and the resulting embedding is unique.

When we come to an edge we have reached before, skip it and move to the next unused edge on that vertex. 
If there are no more unused edges on the vertex, then backtrack to the last vertex with an unused edge. 
Continue in this way until all the edges have been drawn (the procedure will terminate, since the graph is
connected).  If we are able to do this without edges crossing, we will have a planar embedding of the split
graph, uniquely determined by the underlying graph and the vertex orientations.  $\Box$

\begin{figure} [h]
    $$\includegraphics{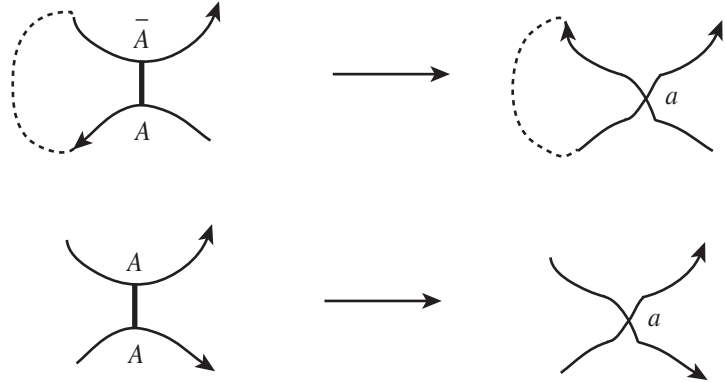}$$
    \caption{Expanding chords into crossings (reversing splits)} \label{F:expand2}
    \end{figure}

Lemma \ref{L:planar} allows us to embed the split graph in the plane.  As in earlier sections, if the original Gauss code is realizable then we can recover an immersion of a graph from this embedding code by connecting the two occurrences of each symbol in the split graph by disjoint chords which respect the orientations of the symbols, and then expanding the chords as in Figure \ref{F:expand2}.  
This leads us to the following theorem.

\begin{thm} \label{T:graph}
Let S be a Gauss code for a connected graph G, with split code $S^*$ and split graph $G^*$.  Then $S$ is realizable if and only if:
\begin{itemize}
   \item $G^*$ is planar.
   \item The embedding of $G^*$ determined by $S^*$ (as in Lemma \ref{L:planar}) can be extended to an embedding of $S^*$.
\end{itemize}
\end{thm}
{\sc Proof:}  If $S$ is realizable, then the split code $S^*$ can be embedded in the plane as an oriented chord diagram.  So the underlying split graph can be embedded (and is therefore planar), and this embedding extended to the chord diagram.  Conversely, if the chord diagram can be embedded, we can expand each chord into a crossing as in Figure \ref{F:expand2}.  It only remains to check that the result is an immersion of $G$ with Gauss code $S$.

As in Theorem \ref{T:closedcurve}, we expand the chords between symbols in the reverse order that the symbols were split.  Say that $S^*$ was produced from $S$ via a series of sequences $S = S_0, S_1, S_2, ..., S_n, S_{n+1} = S^*$, where $S_{i+1}$ is the result of splitting crossing $p_i$ in $S_i$.  Since each splitting was chosen to keep the graph connected, expanding the splits in the reverse order will also keep the graph connected.  Our goal is for these expansions to exactly reverse the effects of the splitting process on the original Gauss code.  There is one possible ambiguity:  the expansion process shown in Figure \ref{F:expand2} assumes that crossings between different edges were always split so as to preserve the orientations of the edges; but, in order to keep the graphs connected, it may have been necessary to reverse the orientations of one of the two edges, as shown in Figure \ref{F:split2}.  So the expansion may leave one of the edges with the ``wrong" orientation.  However, since we have recorded all the sequences $S_i$, this is easily fixed - simply compare the code for our current graph with the corresponding $S_i$.  They will differ in at most the orientation of one edge, so we can simply reverse the orientation of this edge in our graph, and continue the process.  The final result will be an immersion of $G$ with Gauss code $S$.  $\Box$

\begin{cor} \label{C:gauss}
Given a Gauss code $S$ on a connected graph $G$, there is a polynomial time algorithm to determine whether $S$ is realizable.
\end{cor}
{\sc Proof:}  Producing the split code can be done in linear time, and there are linear-time algorithms for determining whether the split graph is planar \cite{bm}.  Finally, by Corollary \ref{C:cut}, we can determine whether the embedding of $G^*$ extends to $S^*$ in polynomial time as well. $\Box$\\

\noindent{\sc Remark:}  In the theory of spatial graphs, a spatial embedding of a graph can be described by the {\it Gauss code} for a particular diagram of the knot (a projection of the knot to the plane).  In this context, the Gauss code comes with additional information at each crossing, indicating which edge crosses over the other and the sign of the crossing.  Corollary \ref{C:gauss} implies that there is a polynomial time algorithm to determine whether a Gauss code is realizable in this context as well.  Once we have determined whether the code can be realized by an immersion of the graph in the plane, and constructed such an immersion, using Theorem \ref{T:graph}, we simply need to check whether the sign and over/under information at each crossing are compatible.

\small

\normalsize

\end{document}